\documentclass{amsart}
\usepackage{amsmath,amscd,amsfonts,amssymb}

\newtheorem{theorem}{Theorem}[section]
\newtheorem{lemma}[theorem]{Lemma}

\newtheorem{corollary}[theorem]{Corollary}
\theoremstyle{definition}

\theoremstyle{remark}

\numberwithin{equation}{section}

\begin{document}



\title[Bounds for the relative n--th nilpotency
degree...]{Bounds for the relative n--th nilpotency degree in
compact groups}

\author[R. Rezaei]{Rashid Rezaei}

\address{Department of Mathematics, University of Malayer, P. O. Box  657719, 95863, Malayer, Iran}

\email{ras$\_$rezaei@yahoo.com}

\author[F.G. Russo]{Francesco G. Russo}

\address{Structural  Geotechnical Dynamics Laboratory StreGa, University of Molise, Via Duca degli Abruzzi, 86039, Termoli (CB), Italy
\newline
and
\newline
Department of Mathematics, University of Palermo, Via Archirafi 34,
90123, Palermo, Italy.}

\email{francescog.russo@yahoo.com}




\date{Received: xxxxxx; Revised: yyyyyy; Accepted: zzzzzz.
\newline \indent $^{*}$ Corresponding author}

\begin{abstract} The line of investigation of the present paper goes back to a
classical work of W. H. Gustafson of the 1973, in which it is
described the probability that two randomly chosen group elements
commute. In the same work,  he gave some bounds for this kind of
probability, providing information on the group structure. We have
recently obtained some generalizations of his results for finite
groups. Here we improve them in the context of the compact
groups.\end{abstract}

\maketitle

\section{Introduction}
A compact group $G$ admits a unique left Haar measure $\mu_G$ which
is normalized and left-invariant (see \cite[Sections 18.1, 18.2,
Proposition 18.2.1]{hr}).  This allows us to assume that $G$ has a
unique probability measure space with respect to $\mu_G$ (see
\cite[Sections 18.1, 18.2]{hr} or \cite[Section 2]{g}). On the
product measure space $G \times G$, it is possible to consider the
product measure $\mu_G\times \mu_G$ which is a probability measure.
If \[C_2=\{(x,y)\in G \times G  \ | \ [x,y]=1\},\] then
$C_2=f^{-1}(1)$, where \[f: (x,y)\in G \times G \mapsto
f(x,y)=[x,y]\in G.\] Clearly, $f$ is continuous and $C_2$ is a
compact measurable subset of $G \times G $. Therefore it is possible
to define \[d(G)=(\mu_G\times \mu_G)(C_2)\] as the $commutativity$
$degree$ of $G$.  In the finite case $d(G)$ is described in
\cite{mog, das-nath1, das-nath2, elr, le, l, les}. We may extend the
notion of $d(G)$ as follows. Suppose that $n\geq1$, $G^n$ is the
product of $n$-copies of $G$ and $\mu^n_G$ that of $n$-copies of
$\mu_G$. We define
\[d^{(n)}(G)=\mu^{n+1}_G(C_{n+1})\] as the $n$-$th$ $nilpotency$
$degree$ of $G$, where
\[C_{n+1}=\{(x_1,\ldots,x_{n+1})\in G^{n+1} \ | \
[x_1,x_2,...,x_{n+1}]=1  \}.\] Obviously, if $G$ is finite, then $G$
is a compact group with the discrete topology and so the Haar
measure of $G$ is the counting measure. Then, for a finite group
$G$, we have
\[d^{(n)}(G)=\mu^{n+1}_G(C_{n+1})= \frac{|C_{n+1}|}{|G|^{n+1}}.\]
See for details \cite{elr, l}.

More generally, let $H$ be a closed subgroup of a compact group $G$.
We define \[D_2=\{(h,g)\in H \times G  \ | \ [h,g]=1\}\] and note
that $D_2=\phi^{-1}(1)$, where \[\phi:(h,g)\in H \times G \mapsto
\phi(h,g)=[h,g]\in G.\] Clearly, $\phi$ is continuous and $D_2$ is a
compact measurable subset of $H \times G $. Note that $\phi$ is the
restriction of $f$ to $H\times G$ and this  shows that $H$ has to be
closed subgroup of $G$, if we want to preserve the topological
structure. Then we define \[d(H,G)= (\mu_H \times \mu_G) (D_2)\] as
the $relative$ $commutativity$ $degree$ of $H$ with respect to $G$.
 Considering
\[D_{n+1}=\{(h_1,, ..., h_n, g)\in H^n \times G  \ |\
[h_1,h_2,...,h_n,g]=1\},\] we define \[d^{(n)}(H,G)= (\mu^n_H \times
\mu_G) (D_{n+1})\] as the $relative$ $n$-$th$ $nilpotency$ $degree$
of $H$ with respect to $G$. As already noted, \cite{mog, elr, le, l,
les} give contributions to the knowledge of the $n$-th nilpotency
degree in case of finite groups. Recently, the case of infinite
groups can be found in \cite{ek, er1, er2, er3, er4, rr, r}. We will
try to extend the results in \cite[Sections 3,4,5]{elr} looking at
the methods in \cite{ek, er1, er2, er3, rr, r}.

\section{Relative commutativity degree}
The next statement is useful for proving most of our results.

\begin{lemma}Assume that $G$ is a compact group, $H$ is a closed subgroup of
$G$ and $C_G([h_1,...,h_n])$ is the centralizer of the commutator
$[h_1, ...,h_n]$ in $G$ for some elements $h_1, ..., h_n$ in $H$.
Then $$ d^{(n)}(H,G)=\int_H \ldots \left( \int_H \mu_G (C_G([h_1,
..., h_n]))d\mu_H (h_1)\right) \ldots
 d\mu_H(h_n),$$ where
\[\mu_G(C_G([h_1,...,h_n]))=\int_G\chi_{_{D_{n+1}}}(h_1,..., h_n,g)
d\mu_G(g)\] and $\chi_{_{D_{n+1}}}$ denotes the characteristic map
of the set $D_{n+1}$.\end{lemma}

\begin{proof} Since \[\mu_G(C_G([h_1,...,h_n]))=\int_G
\chi_{_{D_{n+1}}}(h_1,..., h_n,g)d\mu_G(g),\] we have by
Fubini-Tonelli's Theorem: \[d^{(n)}(H,G)=(\mu^n_H\times
\mu_G)(D_{n+1})= \int_{H^n \times G}
\chi_{_{D_{n+1}}}(d\mu^n_H\times d\mu_G)\] \[=\int_H ... \left(
\int_H \left( \int_G \chi_{_{D_{n+1}}}(h_1,...,h_n,g) d\mu_G(g)
\right) d\mu_H (h_1)\right) ... \ d\mu_H(h_n)\]
\[=\int_H ... \left( \int_H \mu_G(C_G([h_1,...,h_n]) d\mu_H(h_1) \right)... \ d\mu_H(h_n).\]
\end{proof}

We recall the following elementary fact, which can be found in
\cite{g}. See also \cite[Lemma 3.1]{er1}.

\begin{lemma}Assume $H$ is a closed subgroup of a compact group $G$. If $|G
: H| = n< \infty$, then $\mu_G(H)=\frac{1}{n}$. If  $|G :
H|=\infty$, then $\mu_G(H) = 0$.\end{lemma}

\begin{proof} Assume that $|G : H| = n$ is finite. Then $G = {\overset{n}
{\underset{i=1} \bigcup}} g_iH$. So we have \[1 = \mu_G(G) = \mu_G(
\bigcup^n_{i=1} g_iH) = \sum^n_{i=1} \mu_G(g_iH) = \sum^n_{i=1}
\mu_G(H) = n \mu_G(H)\] and therefore $\mu_G(H) = \frac{1}{n}$ . Now
assume that $\alpha=|G : H| = \infty$. Of course, $\alpha>0$, then
$t\alpha > 1$ for some positive integer $t$. By assumption, $G =
{\underset{i\in I} \bigcup} g_iH$, where $I$ is an infinite set.
Choose a subset $J$ of $I$ of cardinality $t$. It follows that
\[1 = \mu_G(G) \geq \mu_G( \bigcup_{j\in J} g_jH) \geq \sum_{j\in J}
\mu_G(g_jH) = t \alpha >0.\] This contradicts $\mu_G(H) = 0$ and the
proof of the lemma follows. \end{proof}

Lemma 2.2 will be used in most of our proofs, even if the following
form is more suitable.

\begin{lemma}Assume $H$ is a closed subgroup of a compact group $G$. If $|G
: H| \geq  n$, then $\mu_G(H)\leq \frac{1}{n}$. If $|G : H| \leq n$,
then $\mu_G(H)\geq \frac{1}{n}$. In particular, $|G : H| = n$ if and
only if $\mu_G(H)= \frac{1}{n}$.\end{lemma}

\begin{proof} This follows from an argument as in Lemma 2.2.
\end{proof}

Lemma 2.3 allows us to reformulate \cite[Theorem 3.10]{elr} for
infinite groups in terms of the following result. The reader may
find exactly the same proof in \cite{rr}: here we repeat it, just
for sake of completeness and because we want to point out the
methods and the ideas which are often used in similar circumstances.

\begin{theorem} Let $H$ be a closed subgroup of a compact group $G$.
\begin{itemize}
\item [(i)] If $d(H,G)=\frac{3}{4}$, then $H/(Z(G)\cap H)$ is cyclic
of order 2.
\item[ii)] If $d(H,G)=\frac{5}{8}$ and $H$ is nonabelian, then $H/(Z(G)\cap H)$ is
2-elementary abelian of rank 2.
\end{itemize}\end{theorem}

\begin{proof} (i). Assume that $d(H,G)=\frac{3}{4}$ and let $K=H \cap Z(G)$.
If $h$ is a element of $H$ not belonging to $K$, then
$|G:C_G(h)|\geq 2$ and so $\mu_G(C_G(h))\leq \frac{1}{2}$ by Remark
2.3. On the other hand, if $h$ is an element of $K$, then
$\mu_G(C_G(h))=1$. From these facts and Lemma 2.1, we have
\[\frac{3}{4}=d(H,G)= \int_H \mu_G (C_G(h))d\mu_H (h)\]
\[=\int_K \mu_G (C_G(h))d\mu_H (h)+ \int_{H-K} \mu_G (C_G(h))d\mu_H (h)\]
\[\leq\int_K d\mu_H (h)+\frac{1}{2}\int_{H-K}d\mu_H(h)=\mu_H(K)+\frac{1}{2}(1-\mu_H(K)).\]
Therefore,  $\mu_H(K)\geq \frac{1}{2}$. On the other hand, $K$ is a
closed subgroup of the abelian group $H$ such that $\mu_H(K)\leq
\frac{1}{2}$. Then $\mu_H(K)= \frac{1}{2}$ and so $|H:K|=2$. This
means that $H/K$ is cyclic of order 2, as claimed.

(ii). Assume that $d(H,G)=\frac{5}{8}$ and let $K=H \cap Z(G)$. We
may argue as in the previous statement (i). On a hand, we have
$\frac{5}{8}=d(H,G)\leq\frac{1}{2}+\frac{1}{2}\mu_H(K).$ Therefore,
$\mu_H(K)\geq \frac{1}{4}$. On the other hand, $K$ is a closed
subgroup of the nonabelian group $H$ so that $\mu_H(K)\leq
\frac{1}{4}$, still by Lemma 2.3 . This gives $\mu_H(K)=
\frac{1}{4}$ so that $|H:K|=4$. This means that $H/K$ has order 4.
Since $H$ is nonabelian,  $H/K$ cannot be cyclic. From this, $H/K$
is 2-elementary abelian of rank 2, as claimed. \end{proof}

Note that \cite[Theorem 3.10]{elr} follows from Theorem 2.4 when we
consider a finite group with the counting measure on it. Now we
extend \cite[Lemma 3.2]{elr} to the case of infinite groups. The
next result overlaps \cite[Lemma 3.2]{er1}.

\begin{lemma}Let $H$ be a closed subgroup of a compact group $G$. Then
\[\mu_G(C_G(x))\leq \mu_H(C_H(x))\] for all $x\in G$.\end{lemma}

\begin{proof}Consider the map \[f : hC_H(x)\in \{hC_H(x) \ | \ h\in H
\}\mapsto f(hC_H(x))=hC_G(x) \in \{gC_G(x) \ | \ g\in G \}.\] $f$ is
one--to--one and so $|H:C_H(x)|\leq |G:C_G(x)|$. This implies
$\mu_G(C_G(x))\leq \mu_H(C_H(x))$. \end{proof} An important
dominance condition is the following.

\begin{theorem}Let $H$ be a closed subgroup of a compact group $G$. Then
\[d(G)\leq d(H,G) \leq d(H).\]\end{theorem}

\begin{proof} From Lemma 2.5, $\mu_H(C_H(x))\geq \mu_G(C_G(x))$. Integrating
over $H$ and keeping in mind Lemma 2.1, we have
\[d(H)=\int_H \mu_H(C_H(x)) d\mu_H (x)\geq \int_H \mu_G
(C_G(x))d\mu_H (x)=d(H,G).\] On the other hand, Lemmas 2.1 and 2.5
give
\[d(H,G)=\int_G\mu_H(C_H(x))d\mu_G(x)\geq\int_G\mu_G(C_G(x))d\mu_G(x)=d(G).\]
\end{proof}

\begin{theorem}Let $H$ be a closed subgroup of a compact group $G$. Then
\begin{itemize}
\item [{\rm(i)}]$d(G)\leq \frac{1}{2}+\frac{1}{2}\mu_G(Z(G))$;
\item[{\rm(ii)}]$d(H,G)\leq \frac{1}{2}+\frac{1}{2}\mu_H(K)$, where $K= H \cap Z(G)$.
\end{itemize}\end{theorem}

\begin{proof} (i). By Lemma 2.1 and noting  that $\mu_G(C_G(x))\leq
\frac{1}{2}$ for each noncentral element $x$ of $G$, we have
\[d(G)=\int_G \mu_G(C_G(x)) d\mu_G (x)\]\[ =\int_{Z(G)} \mu_G(C_G(x))
d\mu_G (x)+\int_{G-Z(G)} \mu_G(C_G(x)) d\mu_G (x)\]\[ =\mu_G
(Z(G))+\int_{G-Z(G)} \mu_G(C_G(x)) d\mu_G (x)\]\[ \leq \mu_G
(Z(G))+\frac{1}{2}
(1-\mu_G(Z(G)))=\frac{1}{2}+\frac{1}{2}\mu_G(Z(G)).\]

(ii). By Lemma 2.1 and noting that $\mu_G(C_G(h))\leq \frac{1}{2}$
for each  element $h$ of $H-K$,
\[d(H,G)=\int_H \mu_G(C_G(h)) d\mu_H (h)\]\[ =\int_{K} \mu_G(C_G(h))
d\mu_H (x)+\int_{H-K} \mu_G(C_G(h)) d\mu_H (h)\]\[
=\mu_H(K)+\int_{H-K} \mu_G(C_G(h)) d\mu_H (h)\]\[\leq \mu_H
(K)+\frac{1}{2}
(1-\mu_H(K))=\frac{1}{2}+\frac{1}{2}\mu_H(K).\]\end{proof}

Note that the upper bounds in \cite[Theorem 3.5]{elr} follow from
Theorem 2.7 when we consider a finite group with the counting
measure on it. The lower bounds in \cite[Theorem 3.5]{elr} cannot be
true in the infinite case, as the infinite dihedral group shows.

\begin{corollary}Assume that $H$ is a closed subgroup of a nonabelian compact
group $G$.
\begin{itemize}
\item [i)]If $H\leq Z(G)$, then $d(H,G)=1.$
\item[(ii)]If $H\not\leq Z(G)$ and $H$ is abelian, then $d(H,G)\leq \frac{3}{4}.$
\item[(iii)]If $H\not\leq Z(G)$ and $H$ is nonabelian, then $d(H,G)\leq \frac{5}{8}.$
\end{itemize}
\end{corollary}

\begin{proof} (i). Obvious.

(ii). Since $H\not\leq Z(G)$, $K=H\cap Z(G)\not\leq H$. As in the
proof of Theorem 2.4 (ii), we have $\mu_H(K)\leq \frac{1}{4}$.
Theorem 2.5 (ii) implies $d(H,G)\leq
\frac{1}{2}+\frac{1}{2}(\frac{1}{4})=\frac{3}{4}$.

(iii). We know from Theorem 2.6 and \cite{g} that $d(H,G)\leq
d(H)\leq \frac{5}{8}$. \end{proof}

Note that \cite[Theorem 3.6]{elr} follows from Corollary 2.8 when we
consider a finite group with the counting measure on it.

\begin{corollary}Let $A$ and $B$ be two closed subgroups of a compact group $G$
such that $A\leq B$. Then $d(A,B)\geq d(A,G) \geq
d(B,G).$\end{corollary}

\begin{proof} As in the proof of Theorem 2.6, the condition
\[|A:C_A(x)|\leq |B:C_B(x)| \leq |G:C_G(x)|\]  implies the condition
$\mu_A(C_A(x))\geq \mu_B(C_B(x))\geq \mu_G(C_G(x))$ for every
element $x$ of $G$. Integrating and keeping in mind Lemma 2.1, we
have \[d(A,B)=\int_A \mu_B (C_B(x))d\mu_A (x)\geq\]
\[d(A,G)=\int_A \mu_G (C_G(x))d\mu_A (x)\geq \int_B \mu_G
(C_G(x))d\mu_B (x)=d(B,G).\]\end{proof}

Note that \cite[Theorem 3.7]{elr} follows from Corollary 2.8 when we
consider a finite group with the counting measure on it. We recall
to convenience of the reader \cite[Lemma 3.8]{elr}.

\begin{lemma} Let $H$ and $N$ be two closed subgroups of $G$ such that $N\leq
H$ and $N$ is normal in $G$. Then $C_H(x)N/N\leq C_{H/N}(xN)$ for
every element $x$ of $G$. Moreover, the equality holds if $N\cap
[H,G]$ is trivial.\end{lemma}

Then we may formulate another interesting dominance condition as
follows.

\begin{theorem}Let $H$ and $N$ be two closed subgroups of a compact group $G$
such that $N\leq H$ and $N$ is normal in $G$. Then $d(H,G)\leq
d(H/N, G/N)d(N).$ In particular, the equality holds if $N\cap [H,G]$
is trivial.\end{theorem}

\begin{proof}Consider $S=\{g\in G \ | \ |H:C_H(g)| \textrm{ is finite}\}$.
We have
\[d(H,G)=\int_G\mu_H(C_H(g))d\mu_G(g)=\int_S\mu_H(C_H(g))d\mu_G(g)\]
\[=\int_S \frac{\mu_H(C_H(g)N)}{|C_H(g)N:C_H(g)|}d\mu_G(g)=\int_S \mu_H(C_H(g)N) \mu_N (C_N(g)) d\mu_G(g).\]
In the last equality we have used the argument just before Theorem
2.4 and the fact that $|C_H(g)N:C_H(g)|$ is finite, getting
\[|C_H(g)N:C_H(g)|=|N:C_H(g)\cap N|=\frac{1}{\mu_N (C_N(g))}.\]
Now we get:
\[d(H,G) \leq \int_G\mu_H(C_H(g)N)\mu_N(C_N(g))d\mu_G(g)\]\[= \int_{\frac{G}{N}}\left(\int_{N}\mu_H(C_H(gx)N)\mu_N(C_N(gx))d\mu_N(x)\right)d\mu_{G/N}(gN).\]
By Lemma 2.10,
$\mu_H(C_H(gx)N)=\mu_{\frac{H}{N}}\Big(\frac{C_H(gx)N}{N}\Big)\leq
\mu_{\frac{H}{N}}(C_{\frac{H}{N}}(gN))$, then
\[d(H,G)\leq \int_{\frac{G}{N}} \mu_{G/N} (C_{G/N}(gN))\left(\int_N \mu_N (C_N
(gx))d\mu_N(x)\right)d\mu_{G/N}(gN).\] On another hand,
\[C_2=\{(x,y)\in N\times N \ | \ [gx,y]=1\}=\{(x,y)\in N\times N \ |
\ gx\in C_G(y)\cap gN\}.\] If $x_0\in C_G(y)\cap gN\neq\varnothing$,
then either $gN=g_0N$ or $g=g_0t$ for some $t\in N$, whence
$C_G(y)\cap gN=g_0(C_G(y)\cap N)=g_0C_N(y)$ and so \[C_2=\{(x,y)\in
N\times N \ | \ gx\in g_0C_N(y)\}=\{(x,y)\in N\times N \ | \ x\in
tC_N(y)\}.\] Therefore
\[\int_{N}\mu_N(C_N(gx))d\mu_N(x)\leq\int_N \mu_N(tC_N(y))d\mu_N(y)\]\[=\int_{N}\mu_N(C_N(y))d\mu_N(y)=d(N).\]
Hence \[d(H,G)\leq
d(N)\int_{\frac{G}{N}}\mu_{\frac{H}{N}}(C_{\frac{H}{N}}(gN))d\mu_{G/N}(gN)=d(N)d(H/N,G/N).\]

In particular, if $N \cap [H,G]=1$, then $C_H(g)=C_H(g)N$ and so
$\mu_H(C_G(g))=\mu_H(C_H(g)N)$ for all  $g\in G.$ Furthermore,  we
have
\[\mu_H(C_H(gn)N)=\mu_{G/N}\Big(\frac{C_H(gn)N}{N}\Big)=\mu_{G/N}(C_{H/N}(gN)).\]
Therefore each inequality becomes equality and so
$d(H,G)=d(H/N,G/N)d(N)$. \end{proof}

\section{Relative n--th commutativity degree}
The present section is devoted to extend some results of Section 2.
For instance, the next statement extends the upper bound in Theorem
2.6.

\begin{theorem}If $H$ is a closed subgroup of a compact group $G$, then
\[d^{(n)}(H,G)\leq d^{(n)}(H).\]\end{theorem}

\begin{proof}We may argue as in the proof of Theorem 2.6 in order to get
\[\mu_H(C_H([h_1, ..., h_n]))\geq \mu_G(C_G([h_1, ..., h_n])).\]
Integrating over $H$ and keeping in mind Lemma 2.1, we have
\[d^{(n)}(H)=\int_H \ldots \left( \int_H \mu_H (C_H([h_1, ...,
h_n]))d\mu_H (h_1)\right) \ldots  d\mu_H(h_n)\]
\[\geq \int_H \ldots \left( \int_H \mu_G (C_G([h_1, ..., h_n]))d\mu_H
(h_1)\right) \ldots  d\mu_H(h_n)=d^{(n)}(H,G).\]
\end{proof}

Note that Theorem 3.1 informs us that the sequence
$\{d^{(n)}(H,G)\}_{n\geq1}$ is increasing for any compact group $G$
and any closed subgroup $H$ of $G$.

The evidences of the finite case and the considerations of many
situations in the infinite case can be summarized in the following
result.

\begin{theorem}If $H$ is a closed subgroup of a compact group $G$ and $K=H\cap
Z(G)$, then $d^{(n+1)}(H,G)\leq
\frac{1}{2}\Big(1+d^{(n)}(H/K)\Big).$\end{theorem}

\begin{proof} Let $A=\{(h_1,...,h_{n+1})\in H^{n+1} \ | \
[h_1,...,h_{n+1}]\in Z(G)\cap H\}$ and $B=H^{n+1}-A$. Then
\[d^{(n+1)}(H,G)=\int_{H^{n+1}}\mu_G([h_1,...,h_{n+1}])d(\mu_H)^{n+1}\]
\[=\int_A\mu_G([h_1,...,h_{n+1}])d(\mu_H)^{n+1}+\int_B\mu_G([h_1,...,h_{n+1}])d(\mu_H)^{n+1}\]
\[\leq\mu_H^{n+1}(A)+\frac{1}{2}\mu_H^{n+1}(B)\leq\mu_H^{n+1}(A)+\frac{1}{2}(1-\mu_H^{n+1}(A))=\frac{1}{2}(1+\mu_H^{n+1}(A)).\] On the other hand,
\[\mu_H^{n+1}(A))=\int_{H}...\int_{H}\mu_{\frac{H}{K}}(C_{\frac{H}{K}}([\bar{h_1},...,\bar{h_n}]))d\mu_H(h_1)...d\mu_H(h_n)\]
\[=\int_{H}...\left(\int_{\frac{H}{K}}\int_K\mu_{\frac{H}{K}}(C_{\frac{H}{K}}([\bar{h_1},...,\bar{h_n}]))d\mu_K(k)d\mu_H(\bar{h_1})\right)...d\mu_H(h_n)\]
\[=\int_{H}...\left(\int_{\frac{H}{K}}\mu_{\frac{H}{K}}(C_{\frac{H}{K}}([\bar{h_1},...,\bar{h_n}]))d\mu_H(\bar{h_1})\right)...d\mu_H(h_n)\]
\[=\int_{\frac{H}{K}}...\int_{\frac{H}{K}}\mu_{\frac{H}{K}}(C_{\frac{H}{K}}([\bar{h_1},...,\bar{h_n}]))d\mu_{\frac{H}{K}}(\bar{h_1})...d\mu_{\frac{H}{K}}(\bar{h_n})
= d^{(n)}(H/K).\] and the result follows.\end{proof}

Note that \cite[Theorem 4.3]{elr} follows from Theorem 3.2 when we
consider a finite group with the counting measure on it. Note that
Theorem 3.2 is true also for groups of the form $A_i\times B_j$,
where $A_i$ is a compact abelian (infinite) group, $B_j$ is a finite
group, $i\in I$ and $j\in J$.

\begin{corollary}If $G$ is a compact group, then
$d^{(n+1)}(G)\leq \frac{1}{2} (1+d^{(n)}(G/Z(G))).$\end{corollary}

\begin{proof} This follows from Theorem 3.2 with $H=G$.
\end{proof}

It is possible to bound $d^{n+1}(G)$  as follows.

\begin{theorem}If $G$ is a compact group, then
\[d^{(n+1)}(G)\leq \frac{1}{2^n} (2^n-1+d(G/Z_n(G))).\]\end{theorem}

\begin{proof} We may repeat the proof of \cite[Theorem 4.5]{elr}, since we
do not need that $G$ is finite. We should only note that
$Z_n(G)/Z(G)=Z_{n-1}(G/Z(G))$ is a closed subgroup of $G/Z(G)$.
\end{proof}

Note that  \cite[Theorem 4.5]{elr} follows from Theorem 3.4 when we
consider a finite group with the counting measure on it. Furthermore
Theorem 3.4 is true for compact groups of the form $A\times B$,
where $A$ is a compact abelian (infinite) group and $B$ is a finite
group. In such a case \cite[Theorem 4.5]{elr} cannot be applied.

\begin{theorem}Let $H$ and $N$ be two closed subgroups of $G$ such that $N\leq
H$ and $N$ is normal in $G$. Then $d^{(n)}(H,G)\leq d^{(n)}(H/N,
G/N).$ In particular, the equality holds if $N\cap [_nH,G]$ is
trivial.\end{theorem}

\begin{proof} Let $\lambda$, $\mu$ and $\nu$ be corresponding Haar measures
on $N$, $G$ and $G/N$ respectively. Consider $S=\{(h_1,...,h_n) \ |
\ |G:C_G([h_1,...,h_n])| \textrm{ is finite}\}$. Then
\[d^{(n)}(H,G)=\int_{H^n}\mu_G(C_G([h_1,...,h_{n}]))d\mu^n_H
=\int_{S}\mu_G(C_G([h_1,...,h_{n}]))d\mu^n_H\]\[=\int_{S}\frac{\mu_G(C_G([h_1,...,h_{n}])N)}{|C_G([h_1,...,h_{n}])N:C_G([h_1,...,h_{n}])|}d\mu^n_H\]
\[=\int_{S}\mu_G(C_G([h_1,...,h_{n}])N)\mu_N(C_G([h_1,...,h_{n}]))d\mu^n_H\]\[\leq\int_{H^n}\mu_G(C_G([h_1,...,h_{n}])N)\mu_N(C_G([h_1,...,h_{n}]))d\mu^n_H\]
\[=\int_{\frac{H}{N}}\int_{N}...\int_{\frac{H}{N}}\int_{N}\mu_G(C_G([h_1a_1,...,h_{n}a_n])N)\mu_N(C_G([h_1a_1,...,h_{n}a_n]))\]
\[d\mu_N(a_1)d\mu_{\frac{H}{N}}(h_1N)...d\mu_N(a_n)d\mu_{\frac{H}{N}}(h_nN).\]
On the other hand,
\[\mu_G(C_G([h_1a_1,...,h_{n}a_n])N)=\mu_{\frac{G}{N}}(\frac{C_G([h_1a_1,...,h_{n}a_n])N}{N})\]\[\leq\mu_{\frac{G}{N}}(C_{\frac{
G}{N}}([h_1N,...,h_{n}N])).\] Therefore
\[d^{(n)}(H,G)\leq\int_{\frac{H}{N}}...\int_{\frac{H}{N}}\mu_{G/N}(C_{\frac{G}{N}}([h_1N,...,h_{n}N]))\]
\[\int_{N}...\int_{N}\mu_N(C_G([h_1a_1,...,h_{n}a_n]))
d\mu_N(a_1)d...d\mu_N(a_n)\mu_{\frac{H}{N}}(h_1N)...d\mu_{\frac{H}{N}}(h_nN)\]
\[\leq\int_{\frac{H}{N}}...\int_{\frac{H}{N}}\mu_{\frac{H}{N}}(C_{\frac{G}{N}}([h_1N,...,h_{n}N]))\mu_{\frac{H}{N}}(h_1N)...d\mu_{\frac{H}{N}}(h_nN)
= d^{(n)}(\frac{H}{N},\frac{G}{N})\]
\end{proof}

An easy consequence is the following.

\begin{corollary}If $N$ is a closed normal subgroup
of a compact group $G$, then $d^{(n)}(G)\leq
d^{(n)}(G/N).$\end{corollary}

\begin{proof} This follows from Theorem 3.5 with $H=G$. \end{proof}

\section{Weakening some bounds}
In the present section we will give some upper and lower bounds for
$d^{(n)}(G)$ and $d^{(n)}(H,G)$ by means of the results which have
been previously found.

\begin{corollary}If G is a compact group which is not nilpotent of class at most
$n$, then $d^{(n)}(G)\leq \frac{2^{n+2}-3}{2^{n+2}}.$\end{corollary}

\begin{proof} $G/Z_{n-1}(G)$ is a nonabelian group by the
 assumptions. From \cite{g}, \[d(G/Z_{n-1}(G))\leq
 \frac{5}{8}.\] Now Theorem 3.4 gives $d^{(n)}(G)\leq \frac{1}{2^{n-1}}
(2^{n-1}-1+\frac{5}{8})=\frac{2^{n+2}-3}{2^{n+2}}. $\end{proof}

 Note that  \cite[Theorem 5.1]{elr} follows from Corollary 4.1 when we
consider a finite group with the counting measure on it.

\begin{corollary}If G is a nontrivial compact group with trivial center, then
$d^{(n)}(G)\leq \frac{2^n-1}{2^n}.$\end{corollary}

\begin{proof}Of course, $Z_n(G)$ is trivial for each $n\geq1$. Thus $G$ is a
nonnilpotent group. In particular $\mu_G(Z(G))=0$ so Theorem 2.7 (i)
implies $d(G)\leq \frac{1}{2}.$ Now the result follows by Theorem
3.4 by induction on $n$. \end{proof}

Our last result improves Corollary 2.8 and extends \cite[Theorem
5.5]{elr}.

\begin{theorem}Assume that $H$ is a proper closed subgroup of a nonabelian
compact group $G$ such that $n\geq1$ and $K=Z(G)\cap H$.
\begin{itemize}
\item [(i)]If $H\leq Z_n(G)$, then $d^{(n)}(H,G)=1.$
\item[(ii)]If $H\not\leq Z_n(G)$ and $H/K$ is a nilpotent group of class at most $n$-$1$, then
$d^{(n)}(H,G)=1.$
\item[(iii)]If $H\not\leq Z(G)$ and $H/K$ is a nonnilpotent group of class at most $n$-$1$, then \[d^{(n)}(H,G)\leq
\frac{2^{n+2}-3}{2^{n+2}}.\]
\end{itemize}\end{theorem}

\begin{proof}(i). This is obvious.

(ii). Of course, if $H/K$ is nilpotent of class at most $n$-$1$,
then $[\bar{x}_1, \ldots, \bar{x}_n]=K$ for some elements
$\bar{x}_1, \ldots, \bar{x}_n$ of $H/K$. Therefore,
$G=C_G([\bar{x}_1, \ldots, \bar{x}_n])$ and we may argue as in
Theorem 3.2, getting $d^{(n)}(H,G)=1.$

(iii). Using Corollary 4.1 and the fact that $H/K$ is nonnilpotent
of class at most $n$-$1$, we get $d^{(n-1)}(H/K)\leq
\frac{2^{n+1}-3}{2^{n+1}}$. By Theorem 3.2, we get
\[d^{(n)}(H,G)\leq \frac{1}{2} (1+d^{(n-1)}(H/K))\leq \frac{1}{2}
(1+\frac{2^{n+1}-3}{2^{n+1}})=\frac{2^{n+2}-3}{2^{n+2}}. \]
\end{proof}



\bibliographystyle{amsplain}

\end{document}